\documentclass[11pt]{article}
\usepackage{graphicx,epsfig}
\usepackage{amsmath,amssymb,euscript}
\usepackage{latexsym}
\usepackage{color}

\newtheorem{theorem}{\bf Theorem}[section]

\newtheorem{proposition}[theorem]{\bf Proposition}

\newtheorem{corollary}[theorem]{\bf Corollary}

\def\NN{{\mathbb N}}

\def\PP{{\mathbb P}}
\def\C{{\mathcal{C}}}
\def\O{{\mathcal{O}}}

\begin{document}
\title{On the filtered polynomial interpolation  at Chebyshev nodes}


\author{Donatella Occorsio{ \thanks{ corresponding author}} and Woula Themistoclakis}

%

%

\maketitle
\begin{abstract}
The paper deals with a special filtered approximation method, which originates  interpolation polynomials at Chebyshev zeros by using de la Vall\'ee Poussin filters. These polynomials can be an useful device for many theoretical and applicative problems since they combine the advantages of the classical Lagrange interpolation, with the uniform convergence in spaces of locally continuous functions equipped with suitable, Jacobi--weighted, uniform norms. The uniform boundedness of the related Lebesgue constants, which equals to the uniform convergence and is missing from Lagrange interpolation, has been already proved in literature under different, but only sufficient, assumptions. Here, we state the necessary and sufficient conditions to get it.
These conditions are easy to check since they are simple inequalities on the exponents of the Jacobi weight defining the norm. Moreover, they are necessary and sufficient to get filtered interpolating polynomials with a near best approximation error, which tends to zero as the number $n$ of nodes tends to infinity. In addition, the convergence rate is comparable with the error of best polynomial approximation of degree $n$, hence the approximation order improves with the smoothness of the sought function.  Several numerical experiments are given in order to test the theoretical results, to make a comparison with the Lagrange interpolation at the same nodes and to show how the Gibbs phenomenon can be strongly reduced.
\end{abstract}



\section{Introduction}
Lagrange interpolation at Chebyshev zeros is one of the most renowned discrete approximation tool, widely used in several fields of the applied sciences.
Also in approximation theory such interpolation process has been widely studied (see e.g. \cite{verteSzabook},  \cite{Sze}, \cite{mastromilobook}, \cite{trefethen}  and the references therein). In particular, concerning the weighted uniform approximation, necessary and sufficient conditions are known in order to get the convergence to sufficiently smooth functions (see e. g. \cite{MaRu},\cite{mastromilobook}). Nevertheless, because of the unboundedness of the associated Lebesgue constants,  such convergence is not ensured for functions with a low degree of smoothness, such as the locally continuous functions on $[-1,1]$.
Moreover, in the case of functions a.e. very smooth apart from some isolated singular points, Lagrange interpolation presents the so-called Gibbs phenomenon producing oscillations and overshoots close to the singularities,  maintained in the regular parts too.

A common method to attenuate Gibbs phenomenon consists in using some filter functions, which generally destroy the interpolating nature of the approximation polynomial. In the present paper we are going to investigate a filtered approximation which preserves the interpolation property at Chebyshev zeros and, at the same time, it provides uniformly bounded Lebesgue constants as well as a strong reduction of the Gibbs phenomenon.

This kind of filtered approximation has its origin in the studies of the Belgian mathematician de la Vall\'ee Poussin concerning the trigonometric polynomial approximation \cite{VPou} and for this reason it  is also known as de la Vall\'ee Poussin (briefly VP) interpolation. The extension of de la Vall\'ee Poussin's work to the algebraic polynomial approximation has been investigated by several authors who defined different kinds of VP algebraic polynomials (see e.g. \cite{Sza, Kis, Ver, Fi.Th1, Hor, Ko, ThBa}). Here we consider the VP interpolating polynomials which have been originated by applying the so--called de la Vall\'ee Poussin filter function to the Fourier--Chebyshev partial sums of degree $n+m-1$, with $0<m<n$. This filtered approximation produces delayed arithmetic means of the Fourier sums having degrees that vary from $n-m$ to $n+m-1$. These VP means are then discretized by means of the Gauss--Chebyshev quadrature rule on $n$ nodes obtaining an algebraic VP polynomial, $V_n^mf$, which interpolates the function $f$ at the same nodes of the applied quadrature rule \cite{woula_99, woula_NA}.

In the limiting case $m=0$ this polynomial coincides with the Lagrange interpolation polynomial, which can be also regarded as a discrete approximation of the Fourier partial sum of order $n-1$.  In the other limiting case $m=n$ it constitutes a discrete and interpolating version of the Fej\'er means that, compared with the simple Fourier sums, have a positive kernel and converge w.r.t. a suitable weighted uniform norm. Nevertheless, they do not share the polynomial projection property and are subject to saturation, i.e. their approximation degree does not increase by increasing the smoothness of the function. By taking $0<m<n$, we get VP interpolating polynomials that, in a certain sense, combine the advantages of the Fourier and Fej\'er operators since both the polynomial preserving property, up the degree $n-m$,  and the convergence w.r.t. weighted uniform norms, are ensured under suitable assumptions. For this reason VP approximation can be an useful device in many situations and it has been already used to get several results in approximation theory as well as in the applications \cite{operatoreD, air, woula_quad,occorsiorusso2011,mastroDB2006,debonis2014,occorsiorussoActa}.

In the present paper we focus on the approximation properties of VP interpolation at Chebyshev zeros w.r.t. weighted uniform norms. In view of the polynomial preserving property, this reduces to investigate  the related Lebesgue constants that have to be  uniformly bounded w.r.t. $n$ and $m$ in order to get a stable, near best polynomial approximation. In literature (see e.g. \cite{woula_99, Ca.Th-wave, woula_NA, ThL1, ThBa, lncs_nostro}) different conditions, only sufficient to get this property, can be found. The main result of the present paper consists in stating the necessary and sufficient conditions in order to get the mentioned property (cf. Theorem \ref{th-LC}). These conditions are simple to be checked since they are inequalities involving the exponents of the Jacobi weight we use to define the weighted uniform norm.

The outline of the paper is as follows. Section 2 deals with the fundamental VP polynomials, which constitutes the basis of the filtered interpolation polynomials studied in Section 3. Finally, in Section 4 we compare VP and Lagrange interpolation at the same nodes. Several numerical experiments 
are given in Section 4 as well as in Section 3.
\section{Fundamental VP polynomials}
Let us consider the four Chebyshev weights
\[
w_1(x):=\frac 1{\sqrt{1-x^2}},\quad w_2(x):=\sqrt{1-x^2},\quad
w_3(x):=\sqrt{\frac{1+x}{1-x}},\quad w_4(x):=\sqrt{\frac{1-x}{1+x}},
\]
and let us adopt the notation without subscript, namely $w(x)$, to mean anyone of the previous weights. Moreover, denote by $\{p_n(w,x)\}_n$ the associated system of orthonormal polynomials with positive leading coefficients.

Setting $t:=\arccos x$, for all $n\in\NN$ it is well known that the n--th orthonormal Chebyshev polynomial has the following trigonometric form
\begin{equation}\label{pol}
p_n(w,x)=\left\{\begin{array}{ll}
\sqrt{\frac 2\pi}\ \cos[nt] \quad \left(\frac 1{\sqrt{\pi}}\ \mbox{for}\ n=0\right)&
\mbox{if $w=w_1$}\\[.15in] \displaystyle
\sqrt{\frac 2\pi}\ \frac{\sin[(n+1)t]}{\sin t} \quad & \mbox{if $w=w_2$}\\[.15in]\displaystyle
\frac 1{\sqrt{\pi}}\ \frac{\cos[(2n+1)t/2]}{\cos[t/2]}  \quad & \mbox{if $w=w_3$}\\[.15in]\displaystyle
 \frac 1{\sqrt{\pi}}\ \frac{\sin[(2n+1)t/2]}{\sin[t/2]}\quad & \mbox{if $w=w_4$}
\end{array}\right.
\end{equation}
where it is understood that we consider the  continuous extension in the not defined cases $t=0,\pi$. Moreover, $p_n(w,x)$ has the following arcsine distributed zeros
\begin{equation}\label{zeros}
x_k=\cos t_k, \qquad\mbox{with}\qquad t_k:=t_{n,k}(w)=\left\{\begin{array}{ll}
\frac{(2k-1)\pi}{2n} & \mbox{if $w=w_1$}\\ [.15in]
\frac{k\pi}{n+1} & \mbox{if $w=w_2$}\\ [.15in]
\frac{(2k-1)\pi}{2n+1} & \mbox{if $w=w_3$}\\ [.15in]
\frac{2k\pi}{2n+1} & \mbox{if $w=w_4$}
\end{array}\right.
\qquad k=1,\ldots,n.
\end{equation}

The fundamental Lagrange polynomials interpolating at these zeros are given by
\begin{equation}\label{fund-lag}
l_{n,k}(x):=\frac{K_n(x,x_k)}{K_n(x_k,x_k)}=\lambda_{n,k}\sum_{j=0}^{n-1}p_j(w,x_k)p_j(w,x),
\qquad k=1,\ldots,n
\end{equation}
where $K_n(x,y):=\sum_{j=0}^np_j(w,x)p_j(w,y)$ denotes the well-known Darboux kernel and
\begin{equation}\label{lambda}
\lambda_{n,k}:=\lambda_{n,k}(w)= \frac 1{\displaystyle\sum_{j=0}^{n}p_j^2(w,x_k)}
=\left\{ \begin{array}{cl}
\frac \pi n &\mbox{if}\quad  w=w_1,\\[.5cm]
\frac{\pi}{n+1}\sin^2t_k &\mbox{if}\quad w=w_2,\\ [.5cm]
\frac{4\pi}{2n+1}\cos^2\frac{t_k}2  &\mbox{if}\quad w=w_3,\\ [.5cm]
\frac{4\pi}{2n+1}\sin^2\frac{t_k}2  &\mbox{if}\quad w=w_4,
\end{array}\right.
\end{equation}
are the well--known Christoffel numbers.

In \cite{woula_99} it was first proved that the interpolation property of fundamental Lagrange polynomials (namely $l_{n,k}(x_h)=\delta_{k,h}$,  $\forall h,k=1,\ldots, n$)
can be also achieved by the so--called VP fundamental polynomials, defined for any pair of degree--parameters $m<n$ in $\NN$, as follows
\begin{equation}\label{fi-filter}
\Phi_{n,k}^m(x)=\lambda_{n,k}\sum_{j=0}^{n+m-1}\mu_{n,j}^m p_j(w,x)p_j(w,x_k),
\end{equation}
where the coefficients $\mu_{n,j}^m$ are the values at $j=0,\ldots, n+m-1$, of the de la Vall\'ee Poussin filter plotted below 

\begin{figure}[h]
\begin{center}
\includegraphics[scale=.45]{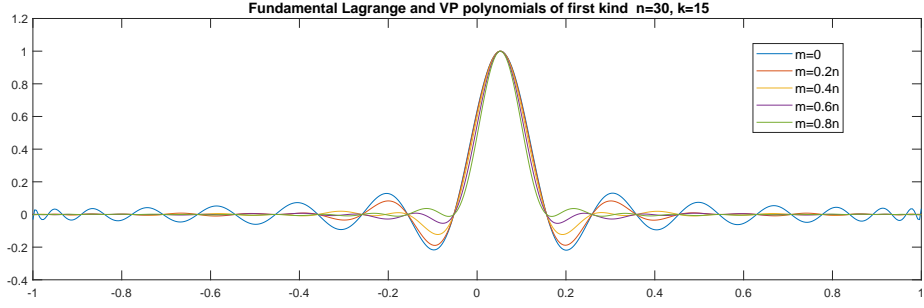}
{\caption{\label{fig-fund}Fundamental Lagrange ($m=0$) and VP polynomials ($m>0$) for $w=w_1$, $n=30$ and $k=15$}}
\end{center}
\end{figure}

\hspace{5cm}
\setlength{\unitlength}{.65cm}
\begin{picture}(7,3.5)
\put(0.3,0.3){\circle*{.125}}\put(0,0){\makebox(0,0){0}}
\put(0.3,0.3){\vector(1,0){7}}\put(7.7,0.3){\makebox(0,0){}}
\put(0.3,0.3){\vector(0,1){2.5}}\put(0,3.3){\makebox(0,0){}}
\put(0,2){\makebox(0,0){1}}\put(0.3,2){\circle*{.125}}
\put(0.3,2){\color{red}\line(1,0){2.5}}
\put(2.8,2){\color{red}\line(1,-1){1.7}}
\put(2.8,0){\makebox(0,0){n-m}}\put(2.8,0.3){\circle*{.125}}
\put(2.8,0.6){\circle*{.08}}\put(2.8,0.9){\circle*{.08}}
\put(2.8,1.2){\circle*{.08}}\put(2.8,1.5){\circle*{.08}}\put(2.8,1.8){\circle*{.08}}
\put(4.5,0){\makebox(0,0){n+m}}\put(4.5,0.3){\circle*{.125}}
\put(4.5,2.2){VP filter}
\end{picture}
\vspace{.5cm}\newline
i.e. defined by
\begin{equation}\label{muj}
\mu_{n,j}^m:=\left\{\begin{array}{ll}
1 & \mbox{if}\quad j=0,\ldots, n-m,\\ [.1in]
\displaystyle\frac{n+m-j}{2m} & \mbox{if}\quad
n-m< j< n+m.
\end{array}\right.
\end{equation}
These filtered polynomials are equivalent to the following discrete VP means of Darboux kernels
\begin{equation}\label{fi-mean}
\Phi_{n,k}^m(x)=\frac {\lambda_{n,k}}{2m}\sum_{j=n-m}^{n+m-1}
K_j(x,x_k),\qquad 0<m<n.
\end{equation}
We also recall the following useful representation of the fundamental VP polynomials \cite{Ca.Th-wave}
\begin{equation}\label{fi-sum}
\Phi_{n,k}^m(x)=\lambda_{n,k}\sum_{j=0}^{n-1} p_j(w,x_k)q_{n,j}^m(w,x),
\end{equation}
where the polynomials
\begin{equation}\label{q-basis}
\hspace{-.2cm}q_{n,j}^m(w,x):=\left\{\begin{array}{ll}
p_j(w,x)&\mbox{if}\quad j=0,\ldots,n\hspace{-.05cm}-\hspace{-.05cm}m,\\ [.1in]
\displaystyle\frac{m\hspace{-.05cm}+\hspace{-.05cm}n\hspace{-.05cm}-
\hspace{-.05cm}j}{2m}p_j(w,x)-\frac{m\hspace{-.05cm}-
\hspace{-.05cm}n\hspace{-.05cm}+\hspace{-.05cm}j}{2m}p_{2n-j}(w,x)
&\mbox{if}\quad n\hspace{-.05cm}-\hspace{-.05cm}m<j<n,
\end{array}\right.
\end{equation}
similarly to $\{p_j(w,x)\}_j$, satisfy the orthogonality relation
\begin{equation}\label{q-ort}
\int_{-1}^1 q_{n,j}^m(w,x)q_{n,i}^m(w,x)w(x)dx=0,\qquad \forall i\ne j\in\{0,\ldots, n-1\} .
\end{equation}
Hence, comparing formulas (\ref{fi-sum}) and (\ref{fund-lag}), we point out that the fundamental VP polynomials are given by the same coefficients of the fundamental Lagrange polynomials, but using the orthogonal polynomial basis $\{q_{n,j}^m\}_{j=0}^{n-1}$ instead of $\{p_j\}_{j=0}^{n-1}$ . Moreover, as we have already mentioned,
these polynomials share the interpolation property (see \cite{woula_99} and \cite{woula_NA})
\begin{equation}\label{int-fund}
\Phi_{n,k}^m(x_h)=l_{n,k}(x_h)=\delta_{k,h}, \qquad \forall h,k=1,\ldots, n,\qquad \forall m<n.
\end{equation}
Indeed, the fundamental VP and Lagrange polynomials coincide in the limiting case $m=0$, while in the case $m>0$ they have different degree, since
\[
\deg\left[\Phi_{n,k}^m\right]=n+m-1, \qquad\mbox{and}\qquad
\deg\left[l_{n,k}\right]=n-1, \qquad\quad\forall k=1,\ldots,n.
\]
We remark that the degree--parameter $n\in\NN$ determines the number of interpolation nodes of both the fundamental polynomials $\Phi_{n,k}^m$ and $l_{n,k}$, while the additional degree--parameter $m<n$ is responsible of increasing the localization of the polynomial  $\Phi_{n,k}^m$ around the node $x_k$. This can be seen in Figure \ref{fig-fund}, which shows how the typical oscillations of the fundamental polynomials are dampened as $m$ grows.

We conclude the section by providing a trigonometric compact formula for each of the four kinds of fundamental VP polynomials.
\begin{proposition}\label{prop}
For any pair of integers $0<m<n$, for each $k=1,\ldots, n$ and for all $x=\cos t$ with $t\in [0,\pi]-\{t_k\ : k=1,\ldots, n\}$ (cf. (\ref{zeros})), we have
\begin{eqnarray}
\label{fi1-trig} w=w_1 &\Longrightarrow&
\Phi_{n,k}^m(x)= \frac{(-1)^{k}}{4 mn}\cos [nt] \Psi(t,t_k),\\
\label{fi2-trig} w=w_2 &\Longrightarrow&
\Phi_{n,k}^m(x)=\frac{(-1)^k\sin t_k}{4 m (n+1)}\ \frac{\sin[(n+1)t]}{\sin t}\ \Psi(t,t_k),\\
\label{fi3-trig} w=w_3 &\Longrightarrow&
\Phi_{n,k}^m(x)=\frac{(-1)^k\cos[t_k/2]}{2 m (2n+1)}\ \frac{\cos[(2n+1)t/2]}{\cos[t/2]}\ \Psi(t,t_k),\\
\label{fi4-trig} w=w_4 &\Longrightarrow&
\Phi_{n,k}^m(x)=\frac{(-1)^k\sin[t_k/2]}{2 m (2n+1)}\
\frac{\sin[(2n+1)t/2]}{\sin[t/2]}\ \Psi(t,t_k),
\end{eqnarray}
where we set
\begin{equation}\label{diff}
\Psi(t,t_k):=\frac{\sin[m(t-t_k)]}{\sin^2[(t-t_k)/2]}
-\frac{\sin[m(t+t_k)]}{\sin^2[(t+t_k)/2]}.
\end{equation}
\end{proposition}
\section{The filtered interpolating polynomials}
Throughout the section we  assume that $f$ is a function  defined everywhere on $]-1,1[$, and suppose $f$ may be also unbounded at the end-points with a known behavior, which is governed by a Jacobi weight
\[
u(x):=v^{\gamma,\delta}(x)=(1-x)^\gamma (1+x)^\delta, \qquad \gamma,\delta \ge 0
\]
such that
\begin{equation}\label{lim}
\lim_{x\rightarrow \pm 1}f(x)u(x)=0,\qquad\mbox{if}\qquad u(\pm 1)=0,\qquad\mbox{respectively.}
\end{equation}
In the case $f$ is also continuous, it is known that it can be uniformly approximated to arbitrary accuracy by polynomials, i.e., denoted by $\PP_n$ the set of the algebraic polynomials of degree at most $n$, we have
\begin{equation}\label{Wei}
\lim_{n\rightarrow\infty}E_n(f)_u=0,\qquad E_n(f)_u:=\inf_{P\in\PP_n}\|(f-P)u\|,
\end{equation}
where we set $\|g\|:=\sup_{|x|\le 1}|g(x)|$.
The speed of convergence to zero of $E_n(f)_u$ depends on the smoothness of $f$ and we refer the reader to the wide existing literature for more details (see e.g. \cite{DT, mastromilobook}).

Here we focus on the polynomials $P$ that can realize a near--best  approximation of $f$, producing an error $\|(f-P)u\|$ comparable with the best approximation order of $E_n(f)_u$.

Supposed to know only the values of $f$ at the Chebyshev nodes (\ref{zeros}), we are going to consider the so--called VP polynomial of $f$, that is defined by means of the fundamental VP polynomials 
introduced in the previous section as follows
\begin{equation}\label{VP}
V_n^mf(x):=\sum_{k=1}^n f(x_k)\Phi_{n,k}^m(x), \qquad |x|\le 1,\qquad m<n.
\end{equation}
Note that $V_n^mf$ is a polynomial of degree at most $n+m-1$, which interpolates $f$ at the zeros of $p_n(w,x)$, since from (\ref{int-fund}) we get
\begin{equation}\label{int-VP}
V_n^mf(x_k)=f(x_k),\qquad \forall k=1,\ldots,n,\quad\mbox{and}\quad \forall m<n.
\end{equation}
In order to estimate the approximation error $\|(f-V_n^mf)u\|$, it is important to recall the following invariance property  (see e.g. \cite[Eq. (17)]{woula_NA})
\begin{equation}\label{inva-VP}
V_n^m P=P,\qquad \quad \forall P\in \PP_{n-m}, \qquad \forall m<n.
\end{equation}
From this polynomial preserving property we easily deduce that
\begin{equation}\label{equi-gen}
E_{n+m-1}(f)_u\le \|(f-V_n^mf)u\|\le \left(1+\|V_n^m\|_u\right) E_{n-m}(f)_u, \qquad \forall m<n,
\end{equation}
where $\|V_n^m\|_u$ denotes the so--called Lebesgue constant
\begin{equation}\label{Leb}
\|V_n^m\|_u:=\sup_{f\ne 0}\frac{\|(V_n^mf)u\|}{\|fu\|}.
\end{equation}
In order to have a stable and near best approximation, the Lebesgue constants have to be uniformly bounded w.r.t. $n$ and $m$, and the following theorem states the necessary and sufficient conditions to get this property. Such a result holds in the case that the degree parameters are such that $m<n\le \C m$ for some fixed constant $\C > 1$ independent of $n,m$. In the sequel, we write $n\sim m$ to denote that $n,m$ are related like this. Moreover, we use $\C$ to denote any positive constant, which may have different values at  different occurrences, and we write $\C\ne\C(m,n,\ldots)$ to mean that $\C>0$ is independent of $n,m, \ldots$.
\begin{theorem}\label{th-LC}
Let $u(x)=(1-x)^\gamma (1+x)^\delta$ be a given Jacobi weight with $\gamma,\delta\ge 0$ and let $w(x)$ be any Chebyshev weight w.r.t. which the VP approximation is considered. We have
\begin{equation}\label{boundedness}
\sup_{n\sim m}\|V_n^m\|_u <\infty,
\end{equation}
if and only if the exponents $\gamma$ and $\delta$ satisfy the following bounds
\begin{equation}\label{eq-LC}
\begin{array}{llllll}
\diamond \ \mbox{Case $w=w_1$:} &  0\le \gamma\le 1 & \mbox{and} & 0\le \delta\le 1 && \\
\diamond\ \mbox{Case $w=w_2$:} & 0< \gamma\le 3/2 & \mbox{and} & 0< \delta\le 3/2 & \mbox{and} & -1\le\gamma-\delta \le 1\\
\diamond\ \mbox{Case $w=w_3$:} & 0\le \gamma\le 1 & \mbox{and} & 0< \delta\le 3/2 & \mbox{and} & \gamma-\delta \le 1/2\\
\diamond\ \mbox{Case $w=w_4$:} & 0< \gamma\le 3/2 & \mbox{and} & 0\le \delta\le 1 & \mbox{and} & \gamma-\delta \ge -1/2
\end{array}
\end{equation}
\end{theorem}
We recall that other (more restrictive) sufficient conditions for (\ref{boundedness}) have been stated in \cite[Th. 2.3]{Ca.Th-wave} and in \cite[Th. 3.2]{woula_NA} while in \cite{ThBa} it has been already proved that the bounds (\ref{eq-LC}) imply (\ref{boundedness}), but using an additional hypothesis if $w=w_1$. Here, Theorem \ref{th-LC} concludes the investigation for all the Chebyshev weights by stating the equivalence between (\ref{eq-LC}) and (\ref{boundedness}). In the particular case $w=w_1$ this result can be also deduced from the bivariate case recently studied in \cite{paperoAMC}.

Theorem \ref{th-LC} gives us a simple way to know the uniform boundedness of the  Lebesgue constants (\ref{Leb}) by checking whether or not the exponents of the Jacobi weight $u$ satisfy the conditions in (\ref{eq-LC}).

\begin{figure}
\centering
\includegraphics[scale=0.4]{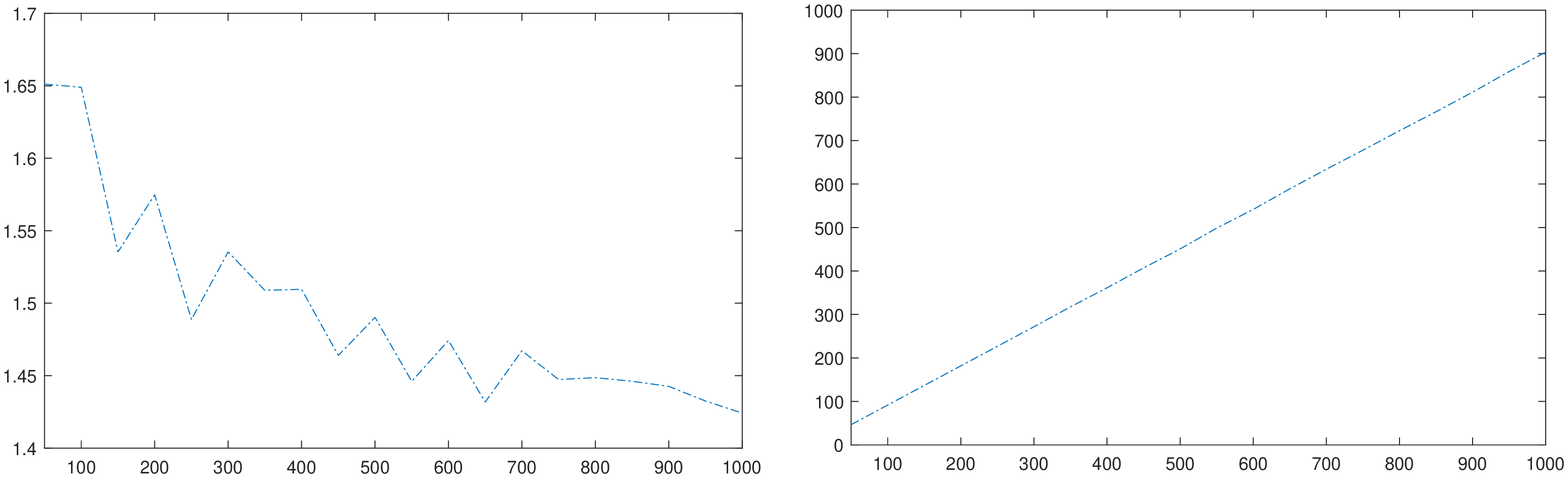}
{\caption{\label{fig1}LC for $w=w_1$, $\theta=0.5$ and $u=v^{\gamma,\delta}$ with  $\gamma=\delta=0.5$ on the left hand side, and  parameters (out of ranges) $\gamma=0.3$, $\delta=1.5$ on the right--hand side}}
\vspace{0.4cm}
\centering
\includegraphics[scale=0.4]{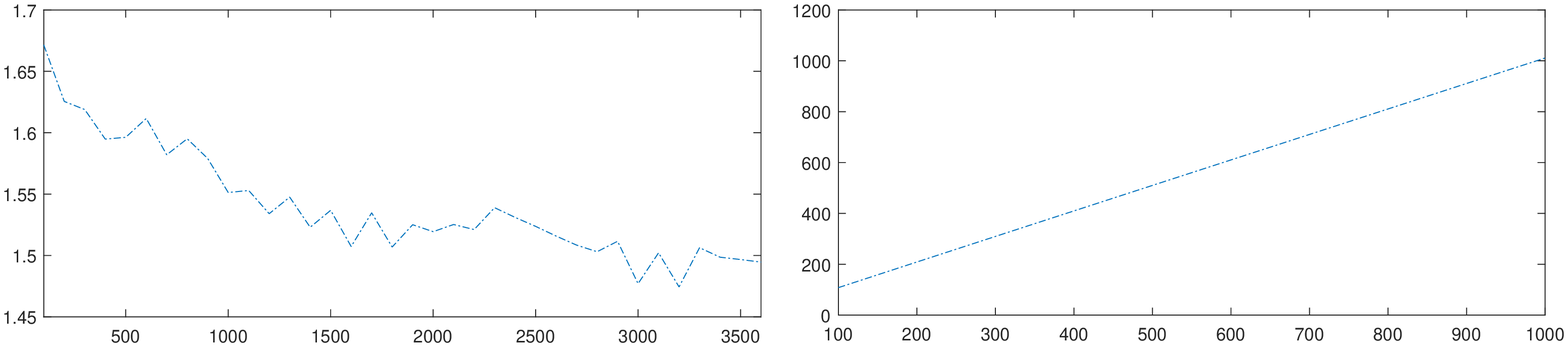}
{\caption{\label{fig2}LC for $w=w_2$, $\theta=0.5$ and $u=v^{\gamma,\delta}$ with  $\gamma=0.7,\ \delta=1.3$ on the left hand side, and  parameters (out of ranges) $\gamma=0$, $\delta=1.5$ on the right--hand side}}
\centering
\includegraphics[scale=0.4]{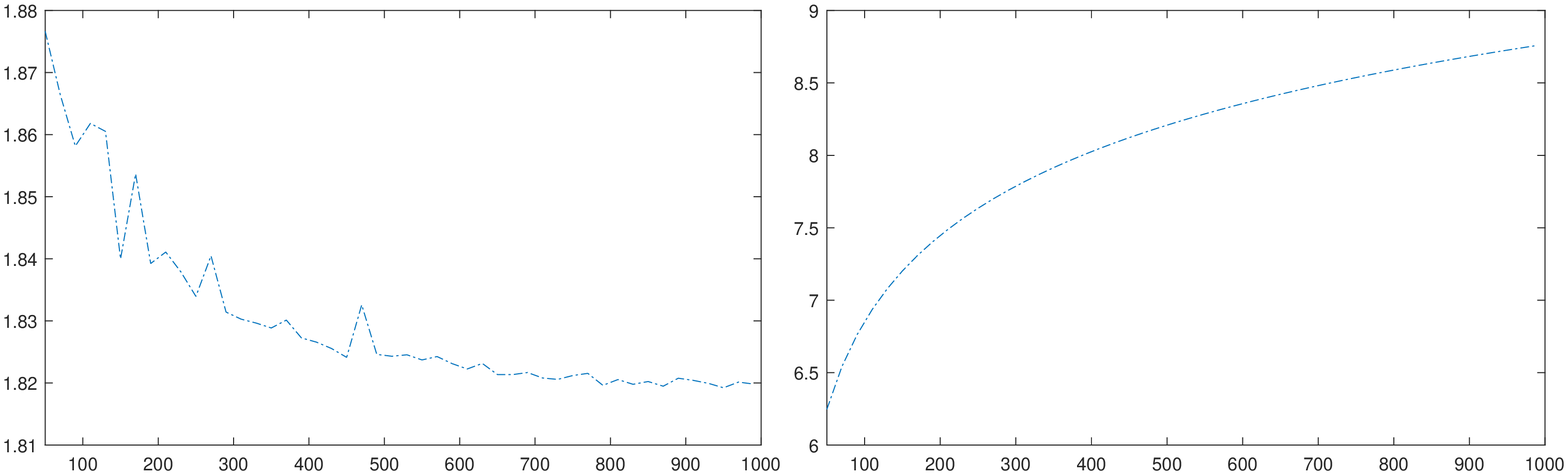}
{\caption{\label{fig3}LC for $w=w_3$, $\theta=0.5$ and $u=v^{\gamma,\delta}$ with  $\gamma=0.4,\ \delta=1.5$ on the left hand side, and  parameters (out of ranges) $\gamma=0.4$, $\delta=0$ on the right--hand side}}
\centering
\includegraphics[scale=0.4]{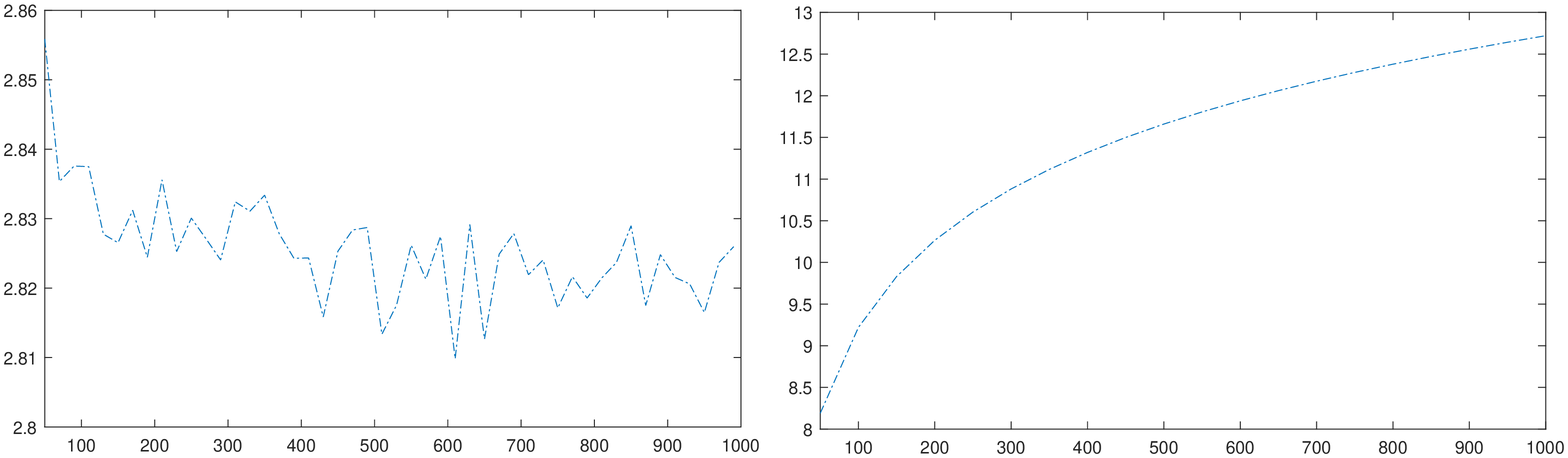}
{\caption{\label{fig4}LC for $w=w_4$, $\theta=0.5$ and $u=v^{\gamma,\delta}$ with  $\gamma=0.5,\ \delta=1$ on the left hand side, and  parameters (out of ranges) $\gamma=0$, $\delta=0.5$ on the right--hand side}}
\end{figure}

In Figures \ref{fig1}--\ref{fig4} we investigated, for each Chebyshev weight, the behavior of the weighted Lebesgue constants (LC) of the interpolating VP polynomial in two cases where the bounds (\ref{eq-LC}) hold (plots on the left hand side) and do not hold (plots on the right hand side) for the exponents of the Jabobi weight defining the weighted uniform norm.
For the numerical computation of the Lebesgue constants (\ref{Leb}) we used the formulas
\begin{equation}\label{Leb-VP1}
\|V_n^m(w)\|_u=\sup_{|x|\le 1}u(x)\sum_{k=1}^n\frac{|\Phi_{n,k}^m(x)|}{u(x_k)},\qquad
\Phi_{n,k}^m(x)=\lambda_{n,k}\sum_{j=0}^{n-1} p_j(w,x_k)q_{n,j}^m(w,x),
\end{equation}
and we considered increasing degree--parameters $n\in\NN$. In order to be sure that $m\sim n$, we have always fixed
\[
0<\theta<1 \qquad \mbox{and}    \qquad m=\lfloor\theta n\rfloor
\]
The plots in Figures \ref{fig1}--\ref{fig4} correspond to $\theta=0.5$ and confirm the theoretical results of Theorem \ref{th-LC}.

Increasing choices of the parameter $\theta\in ]0,1[$ have been investigated in Table \ref{table_new}, where we numerically computed
\[
\sup_{\scriptsize {\begin{array}{c}
n\in\NN\\
m=\lfloor\theta n\rfloor
\end{array}}}\|V_n^m(w)\|_u
\]
for the fourth Chebyshev weights and  several Jacobi weights $u=v^{\gamma,\delta}$ whose exponents satisfy the bounds in (\ref{eq-LC}). As expected, the computations give decreasing values as the parameter $\theta$, and hence $m$, increases.
\vspace{.1cm}\newline
\begin{table}[h]
\centering{
\begin{tabular}{ll||l|l|l|l|l|l|l|l|l}
{} \hspace{1.5cm}$\theta $ & & $0.1$  & $0.2$ & $0.3$ & $0.4$ & $0.5$ & $0.6$  & $0.7$  & $0.8$  & $0.9$ \\  \hline 
&   \\ [-.4cm]
$w=w_1$& $u=v^{0,0}$  & $2.43$  & $1.99$ & $1.73$ & $1.53$ & $1.42$ & $1.26$  & $1.16$  & $1.10$  & $1.04$
\\ [.4cm]
$w=w_2$&$u=v^{0.5,\ 0.4}$  & $2.42$  & $1.98$ & $1.72$ & $1.52$ & $1.41$ & $1.26$  & $1.16$  & $1.10$  & $1.04$
\\ [.4cm]
$w=w_3$&$u=v^{0.2,\ 0.3}$  & $2.47$  & $2.03$ & $1.76$ & $1.55$ & $1.44$ & $1.28$  & $1.18$  & $1.11$  & $1.06$
\\ [.4cm]
$w=w_4$& $u=v^{0.5,\ 0.4}$  & $2.44$  & $2.00$ & $1.74$ & $1.54$ & $1.43$ & $1.28$  & $1.17$  & $1.10$  & $1.05$\vspace{.2cm}
\end{tabular}}
\caption{\label{table_new}{$\sup_{n\in\NN,\ m=\lfloor\theta n\rfloor}\|V_n^m(w)\|_u $ for different choices of $u,w$ and increasing $\theta$}}
\end{table}

We conclude the section by the following corollary which estimates the error of the $VP$ filtered interpolation in terms of the error of best polynomial approximation.
\begin{corollary}\label{cor}
Let $u(x)=(1-x)^\gamma (1+x)^\delta$ be a Jacobi weight whose exponents satisfy the bounds in (\ref{eq-LC}). Then  for any pair of degree--parameters $n\sim m$, we have
\begin{equation}\label{equi-VP}
E_{n+m-1}(f)_u\le \|(f-V_n^mf)u\|\le \C E_{n-m}(f)_u, \qquad \C\ne \C(n,m,f),
\end{equation}
\end{corollary}
This corollary is an immediate consequence of  Theorem \ref{th-LC} that, combined with (\ref{equi-gen}), gives us the constant $\C=1+\|V_n^m(w)\|_u$ in (\ref{equi-VP}). Choosing $ m=\lfloor\theta n\rfloor$, we have to take into account that if $\theta$ increases then $\C$ becomes smaller, but $E_{n-m}(f)_u$ increases since its approximation degree $(1-\theta)n$  decreases w.r.t. $\theta$.
\section{Filtered vs Lagrange interpolation}
In this section we compare the filtered interpolating  polynomial $V_n^mf$ with the well--known Lagrange polynomial interpolating $f$ at the same nodes, namely
\begin{equation}\label{Lag}
L_nf(x):=\sum_{k=1}^n f(x_k)l_{n,k}(x), \qquad |x|\le 1,
\end{equation}
where the fundamental Lagrange polynomials have been defined in (\ref{fund-lag}).

Obviously, by (\ref{int-fund}), we have
\[
L_nf(x_k)=V_n^mf(x_k)=f(x_k),\qquad k=1,\ldots,n,
\]
but the map $V_n^m:f\rightarrow V_n^mf\in \PP_{n+m-1}$ results to be a polynomial quasi--projection which preserves the polynomials up to the degree $n-m$, while $L_n:f\rightarrow L_nf\in\PP_{n-1}$ is a projection on $\PP_{n-1}$, i.e. $L_n P=P$ holds $\forall P\in\PP_{n-1}$. Consequently, by Faber's theorem, $L_n$ has unbounded Lebesgue constants that grow with $n$, in the best cases as $\log n$ (see also \cite{vertesisolo}). In particular it has been proved that (see e. g. \cite{MaRu},\cite{mastromilobook})
\begin{equation}\label{equi-Lag}
E_{n-1}(f)_u\le \|(f-L_nf)u\|\le \C \log n E_{n-1}(f)_u, \qquad \C\ne \C(n,f),
\end{equation}
holds if and only if the exponents of the Jacobi weight $u$ satisfy the following bounds
\begin{equation}\label{cond-Lag}
\begin{array}{clcll}
0\le \gamma\le 1 & \mbox{and} & 0\le \delta\le 1 & \mbox{if}& w=w_1\\
1/2\le \gamma\le 3/2 & \mbox{and} & 1/2\le \delta\le 3/2 & \mbox{if} & w=w_2\\
0\le \gamma\le 1 & \mbox{and} & 1/2\le \delta\le 3/2 & \mbox{if}& w=w_3\\
1/2\le \gamma\le 3/2 & \mbox{and} & 0\le \delta\le 1 & \mbox{if}& w=w_4
\end{array}
\end{equation}
Comparing (\ref{cond-Lag}) and (\ref{equi-Lag}) with (\ref{eq-LC}) and (\ref{equi-VP}), we can observe more restricted bounds  in (\ref{cond-Lag}) and we find the additional $\log n$ factor in (\ref{equi-Lag}), which is not removable. Nevertheless, taking into account that, as $n\rightarrow \infty$,  $\log n$ tends to infinity very slowly, we can say that in the practice, for finite degrees $n$, the maximum weighted approximation error of Lagrange and VP polynomials are almost comparable.

In fact, we have computed the errors
\[
\mathcal{E}^{VP}_{n,m}(f)_u:=\|(f-V_n^mf)u\|\qquad\mbox{and}\qquad \mathcal{E}^{Lag}_{n}(f)_u:=\|(f-L_nf)u\|,
\]
for the following test functions and related weights
\[
\begin{array}{lll}
\displaystyle f_1(x)=\frac {|x+0.5|^\frac 7 2 \sin(x)}{1+x^2}, & w=w_3, & u=v^{0.6,0.6}, \\ [.2in]
\displaystyle
f_2(x)=
\begin{cases}
\frac{32}9(1+x), & -1\le x\le -\frac 13\\
(1-x)^3, & -\frac 13<x\le \frac 13\\
\frac 4{9}(1-x), & \frac 13<x\le 1
 \end{cases}, & w=w_1& u=v^{0,0},\\ [.4in]
\displaystyle f_3(x)=(1+|x|)^\frac 1 4,  & w=w_4 & u=v^{0.5,0.5}.
 \end{array}
 \]
Concerning the error of best uniform weighted approximation of these functions, we have
\[
E_n(f_1)_u= \O (n^{-3}),\qquad E_n(f_2)_u= \O (n^{-1}),\qquad E_n(f_3)_u= \O (n^{-1}),
\]
and the theoretical estimates are coherent with the numerical results displayed in
Tables \ref{table1} and \ref{table22} (first three columns). Here we note that the errors of Lagrange and VP interpolation are almost comparable, especially in the case of sufficiently smooth functions as $f_1$.
\begin{table}[h]
\centering
	\begin{tabular}{|c|c|c||c|c|c|}
            \hline   $n$  & $\mathcal{E}^{VP}_{n,m}(f_1)_u$ & $\mathcal{E}^{Lag}_{n}(f_1)_u$ &
             $n$  & $\mathcal{E}^{VP}_{n,m}(f_2)_u$ & $\mathcal{E}^{Lag}_{n}(f_2)_u$ \\
			\hline  50 &  $3.7e-7$ & $4.1e-7 $ & 400 & 2.9e-2 & 1.9e-1 \\
			\hline  150 &  $9.9e-9$ & $9.9e-9 $ & 1200 & 9.5e-2 & 1.9e-1 \\
			\hline 250 &  $3.6e-9$ & $2.2e-9 $ & 1600 & 3.8e-3 & 1.3e-1  \\
			\hline 350 &  $4.7e-10$ & $4.7e-10 $ & 2800 & 2.0e-3 & 8.9e-2  \\
            \hline  450 &  $2.0e-10$ & $2.1e-10 $ & 3800 & 3.9e-3& 2.7e-2 \\
           \hline
	\end{tabular}
\caption{\label{table1} Errors of the function $f_1(x)$, with $w=w_3$, $u=v^{0.6,0.6}$, $\theta=0.4$,  and the function $f_2$ with $w=w_1$, $u=v^{0,0}$, $\theta=0.9,$}
\end{table}
\begin{table}[h]
\centering
	\begin{tabular}{|c|c|c||c|c|c|}
            \hline
             $n$  & $\mathcal{E}^{VP}_{n,m}(f_3)_u$ & $\mathcal{E}^{Lag}_{n}(f_3)_u$ & $n$  & $\mathcal{E}^{VP}_{n,m}(f_4)_u$ & $\mathcal{E}^{Lag}_{n}(f_4)_u$ \\
			\hline  51 & 3.4e-3 & 3.3e-3 & 100 &  $1.7e-2$ & $4.1e-2 $ \\
			\hline 101 & 1.8e-3  &  1.7e-3 & 600 &  $2.6e-3$ & $2.5e-1 $ \\
			\hline 201 & 9.0e-4  & 8.8e-4 & 1100 &  $1.2e-3$ & $4.7e-1 $ \\
			\hline  301 & 6.0e-4  & 5.9e-4 & 1600 &  $6.7e-4$ & $6.7e-1 $ \\
            \hline  401 & 4.5e-4  & 4.4e-4 & 2100 &  $3.8e-4$ & $8.7e-1$ \\
             \hline 501 & 3.6e-4  & 3.5e-4 & 2600 &  $2.1e-4$ & $1.1e-0 $ \\
              \hline 601 & 2.3e-4  & 2.2e-4 & 3100 &  $1.0e-4$ & $1.3e-0 $ \\
           \hline
	\end{tabular}
\caption{ \label{table22} {Errors of the function $f_3$, for  $w=w_4$,  $u=v^{0.5,0.5}$,  $\theta=0.3$, and the function $f_4(x)=|x|,$ for $w=w_2$, $u=v^{0.1,0.1}$, $\theta=0.9$}}
\end{table}

\begin{figure}[h]
\centering
\includegraphics[width=5.5cm]{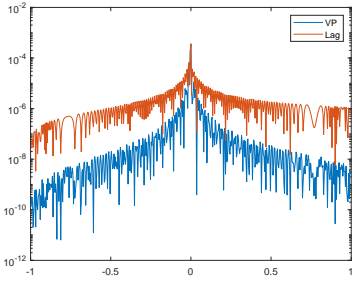}
{\caption{\label{fig5} Pointwise errors by Lagrange (orange)  and VP (blue) for  $f_3$ with $w=w_4$,  $u=v^{0.5,0.5}$,  $\theta=0.3$ }}
\end{figure}

On the other hand,  in the case of the less regular function $f_3$  we plotted the pointwise errors of Lagrange (orange)  and VP (blue) polynomials in  Figure \ref{fig5}, where we can see that locally the
VP polynomial provides an approximation much better than the Lagrange one.

Another point of reflection is based on the comparison between the bounds (\ref{eq-LC}) and (\ref{cond-Lag}). Indeed,  except for the case $w=w_1$, the possible range for the exponents of $u$,  in the case of VP interpolation, is  wider than that one of the corresponding  Lagrange case.
For instance, in Table \ref{table22} (last three columns) we approximated  the function $f_4(x)=|x|$,  with $u=v^{0.1,0.1}$, by the Lagrange and VP interpolating polynomials w.r.t  $w=w_2$. In this case $E_n(f_4)_u=\O(n^{-1})$, but the conditions assuring optimal Lebesgue constants are not satisfied for Lagrange interpolation while they are satisfied for VP interpolation polynomial, where we fixed $m=\theta n$ and $\theta=0.9$. This fact is confirmed by the Lagrange and VP errors we computed in Table \ref{table22} for $f_4$.

\begin{figure}[h]
\centering
\includegraphics[scale=0.5]{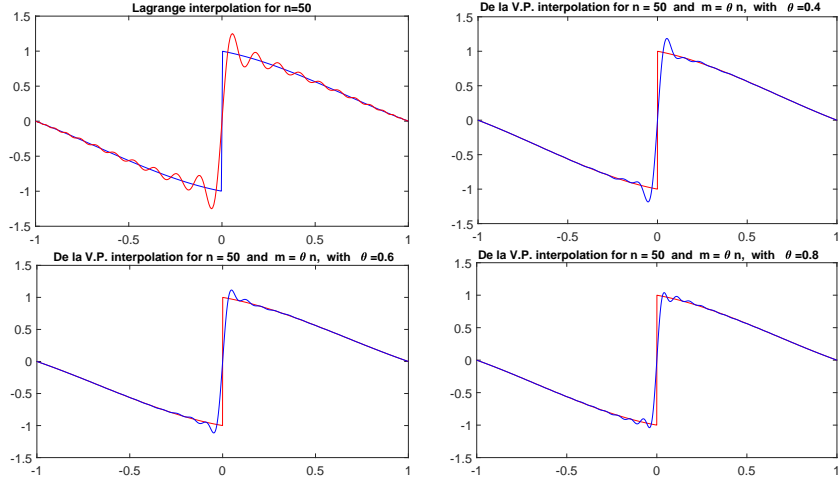}
{\caption{\label{figJUMP} Plots of $uf_5$ (blue line)  and  the weighted approximations $uL_{50}f_5$ and $uV_{50}^{m}f_5$ (in red) for $w=w_2$ , $u=v^{1,1}$ and $m=50 \theta$}}
%
\end{figure}

Finally,  we want to highlight another "good" feature offered by  VP vs Lagrange polynomials when we interpolate bounded variation functions  having  jump discontinuities inside the interval $[-1,1]$. For instance, we consider the function $f_5(x)=\textrm{sign}(x)-\frac x 2$ and the weights $w=w_2$ and $u=v^{1,1}$. In  the first plot of Figure  \ref{figJUMP} we displayed the weighted  function $uf_5$ (blue line) together with the weighted approximations  $uL_{n}f_5$  (red line) for fixed $n=50$. As it was expected,  the Gibbs phenomenon appears and we have  damped oscillations close to the jump, but also along the whole interval. In the other plots of Figure \ref{figJUMP} we considered  $uV_n^mf_5$ of the same degree $n=50$, with several $m=\theta n$ corresponding to the different choices of the parameter $\theta= 0.4$, $0.6$ and $0.8$. As the graphics show, the filtered VP interpolation  induces a  reduction of the Gibbs phenomenon as $\theta$ increases.

Such a reduction of the Gibbs phenomenon can be explained by the major localization provided by the fundamental VP polynomials $\Phi_{n,k}^m$ w.r.t. the fundamental Lagrange polynomials $l_{n,k}$ (cf. Figure \ref{fig-fund}).

\begin{figure}[h]
\centering
\includegraphics[scale=0.5]{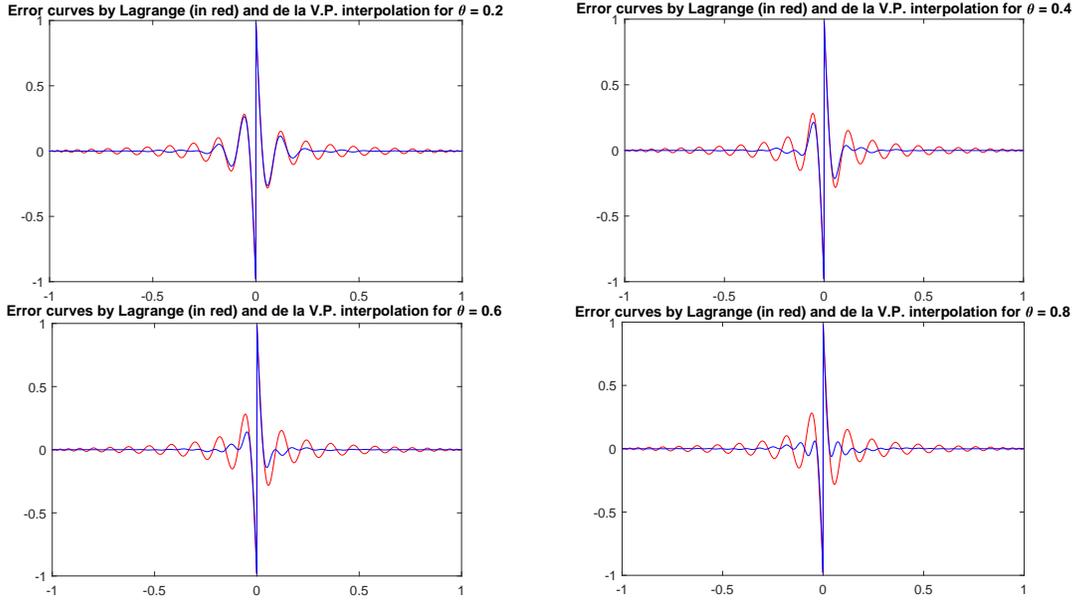}%
{\caption{\label{figERR} Plots of the error curves $u(x)|f_5(x)-L_{50}f_5(x)|$ (red line)  and $u(x)|f_5(x)-V_{50}^mf_5(x)|$ (blue line) for $w_2$ , $u=v^{1,1}$, $m=50 \theta$ and $\theta=0.2$ (up on the left), $\theta=0.4$ (up on the right),  $\theta=0.6$ (down on the left) and $\theta=0.8$ (down on the right)}}
\end{figure}

It is displayed also in Figure \ref{figERR} where the pointwise weighted errors of the proposed VP interpolants (red line), in comparison with the Lagrange ones (blue line) at the same $n=50$ nodes, have been plotted for several $m=\theta n$ corresponding to $\theta=0.2,\ 0.4,\ 0.6,\ 0.8$.

\section{Proofs}
\subsection{Proof of Proposition \ref{prop}}
The statement is based on the following trigonometric identities (see e.g. \cite{Bary})
\begin{eqnarray}
\label{id1}
2\cos\alpha\cos\beta &=&\cos[\alpha-\beta]+\cos[\alpha+\beta],\\ [.1in]
\label{id2}
 2\sin\alpha\sin\beta &=&\cos[\alpha-\beta]-\cos[\alpha+\beta]\\ [.1in]
 \label{id3}
\frac 12+\sum_{j=1}^r\cos [j\alpha]&=&\frac{\sin[(2r+1)\alpha/2]}{2\sin [\alpha/2]},\\
\label{id4}
\sum_{j=0}^r \cos [(2j+1)\alpha/2]&=&\frac{\sin[(r+1)\alpha]}{2\sin [\alpha/2]},\\
\label{id5}
\sum_{r=n-m}^{n+m-1}\sin\left[(2r+1)\alpha/2\right]&=&\frac{\sin [n\alpha] \sin [m\alpha]}{\sin [\alpha/2]},\\
\label{id6}
\sum_{r=n-m}^{n+m-1}\sin\left[(r+1)\alpha\right]&=&\frac{\sin [(2n+1)\alpha/2] \sin [m\alpha]}{\sin [\alpha/2]},
\end{eqnarray}
More precisely, starting from (\ref{fi-mean}) and (\ref{lambda}), we get
\vspace{.1cm}\newline
$\bullet${\it Case $w=w_1$}
\begin{eqnarray*}
\Phi_{n,k}^m(\cos t)&=&
\frac 1 {nm}\sum_{r=n-m}^{n+m-1}\left(\frac 12 +\sum_{j=1}^r\cos[jt]\cos[jt_k]\right)
\\
&\begin{array}{c}\mbox{\footnotesize(\ref{id1})}\\=\end{array}&
\frac 1{2nm}\sum_{r=n-m}^{n+m-1}\left(1 +\sum_{j=1}^r \cos[j(t-t_k)]+\cos[j(t+t_k)]\right)
\\
&\begin{array}{c}\mbox{\footnotesize(\ref{id3})}\\=\end{array}&
\frac 1{4 nm}\sum_{r=n-m}^{n+m-1}\left(\frac{\sin[(2r+1)(t-t_k)/2]}{\sin[(t-t_k)/2]}
+ \frac{\sin[(2r+1)(t+t_k)/2]}{\sin[(t+t_k)/2]}\right)
\\
&\begin{array}{c}\mbox{\footnotesize(\ref{id5})}\\=\end{array}&
\frac 1{4 nm}\left(\frac{\sin[n(t-t_k)]\sin[m(t-t_k)]}{\sin^2[(t-t_k)/2]}
+ \frac{\sin[n(t+t_k)]\sin[m(t+t_k)]}{\sin^2[(t+t_k)/2]}\right)
\end{eqnarray*}
and (\ref{fi1-trig}) follows observing that
\[
\sin[n(t\pm t_k)]=\sin \left[nt \pm (2k-1)\frac{\pi}2\right]=\mp(-1)^{k}\cos[nt]
\]
$\bullet${\it Case $w=w_2$}
\begin{eqnarray*}
&&\Phi_{n,k}^m(\cos t)=
\frac{\sin^2t_k}{m(n+1)}
\sum_{r=n-m}^{n+m-1}\left(\sum_{j=0}^r\frac{\sin[(j+1)t]\sin[(j+1)t_k]}{\sin t\sin t_k}\right)
\\
&\begin{array}{c}\mbox{\footnotesize(\ref{id2})}\\=\end{array}&
\frac{\sin t_k}{2m(n+1)\sin t}
\sum_{r=n-m}^{n+m-1}\left(\sum_{j=1}^{r+1}\cos[j(t-t_k)]-\cos[j(t+t_k)]\right)
\\
&\begin{array}{c}\mbox{\footnotesize(\ref{id3})}\\=\end{array}&
\frac{\sin t_k}{4m(n+1)\sin t}
\sum_{r=n-m}^{n+m-1}\left(\frac{\sin[(2r+3)(t-t_k)/2]}{\sin[(t-t_k)/2]}
- \frac{\sin[(2r+3)(t+t_k)/2]}{\sin[(t+t_k)/2]}\right)
\\
&\begin{array}{c}\mbox{\footnotesize(\ref{id5})}\\=\end{array}&
\frac{\sin t_k}{4m(n+1)\sin t}
\left(\frac{\sin[(n+1)(t-t_k)]\sin[m(t-t_k)]}{\sin^2[(t-t_k)/2]}
- \frac{\sin[(n+1)(t+t_k)]\sin[m(t+t_k)]}{\sin^2[(t+t_k)/2]}\right)
\end{eqnarray*}
and (\ref{fi2-trig}) follows from
\[
\sin[(n+1)(t\pm t_k)]=\sin \left[(n+1)t \pm k\pi\right]=(-1)^{k}\sin[(n+1)t]
\]
$\bullet$ {\it Case $w=w_3$}
\begin{eqnarray*}
&&\Phi_{n,k}^m(\cos t)=
\frac{2\cos^2[t_k/2]}{m(2n+1)}
\sum_{r=n-m}^{n+m-1}\left(\sum_{j=0}^r\frac{\cos[(2j+1)t/2]\cos[(2j+1)t_k/2]}{\cos [t/2]
\cos [t_k/2]}\right)
\\
&\begin{array}{c}\mbox{\footnotesize(\ref{id1})}\\=\end{array}&
\frac{\cos [t_k/2]}{m(2n+1)\cos [t/2]}
\sum_{r=n-m}^{n+m-1}\left(\sum_{j=1}^{r+1}\cos[(2j+1)(t-t_k)/2]+\cos[(2j+1)(t+t_k)/2]\right)
\\
&\begin{array}{c}\mbox{\footnotesize(\ref{id4})}\\=\end{array}&
\frac{\cos [t_k/2]}{2m(2n+1)\cos [t/2]}
\sum_{r=n-m}^{n+m-1}\left(\frac{\sin[(r+1)(t-t_k)]}{\sin[(t-t_k)/2]}
+ \frac{\sin[(r+1)(t+t_k)]}{\sin[(t+t_k)/2]}\right)
\\
&\begin{array}{c}\mbox{\footnotesize(\ref{id6})}\\=\end{array}&
\frac{\cos [t_k/2]}{2m(2n+1)\cos [t/2]}
\left(\frac{\sin[(2n+1)(t-t_k)/2]\sin[m(t-t_k)]}{\sin^2[(t-t_k)/2]} +\right.\\ [.1in]
&& \left.+ \frac{\sin[(2n+1)(t+t_k)/2]\sin[m(t+t_k)]}{\sin^2[(t+t_k)/2]}\right)
\end{eqnarray*}
and (\ref{fi3-trig}) follows from
\[
\sin[(2n+1)(t\pm t_k)/2]=\sin \left[(2n+1)t/2 \pm (2k-1)\frac \pi 2 \right]=\mp(-1)^{k}\cos[(2n+1)t/2]
\]
$\bullet${\it Case $w=w_4$}
\begin{eqnarray*}
&&\Phi_{n,k}^m(\cos t)=
\frac{2\sin^2[t_k/2]}{m(2n+1)}
\sum_{r=n-m}^{n+m-1}\left(\sum_{j=0}^r\frac{\sin[(2j+1)t/2]\sin[(2j+1)t_k/2]}{\sin [t/2]
\sin [t_k/2]}\right)
\\
&\begin{array}{c}\mbox{\footnotesize(\ref{id2})}\\=\end{array}&
\frac{\sin [t_k/2]}{m(2n+1)\sin [t/2]}
\sum_{r=n-m}^{n+m-1}\left(\sum_{j=0}^{r}\cos[(2j+1)(t-t_k)/2]-\cos[(2j+1)(t+t_k)/2]\right)
\\
&\begin{array}{c}\mbox{\footnotesize(\ref{id4})}\\=\end{array}&
\frac{\sin [t_k/2]}{2m(2n+1)\sin [t/2]}
\sum_{r=n-m}^{n+m-1}\left(\frac{\sin[(r+1)(t-t_k)]}{\sin[(t-t_k)/2]}
- \frac{\sin[(r+1)(t+t_k)]}{\sin[(t+t_k)/2]}\right)
\\
&\begin{array}{c}\mbox{\footnotesize(\ref{id6})}\\=\end{array}&
\frac{\sin [t_k/2]}{2m(2n+1)\sin [t/2]}
\left(\frac{\sin[(2n+1)(t-t_k)/2]\sin[m(t-t_k)]}{\sin^2[(t-t_k)/2]} +\right.\\ [.1in]
&& \left.- \frac{\sin[(2n+1)(t+t_k)/2]\sin[m(t+t_k)]}{\sin^2[(t+t_k)/2]}\right)
\end{eqnarray*}
and (\ref{fi4-trig}) follows from
\[
\sin[(2n+1)(t\pm t_k)/2]=\sin \left[(2n+1)t/2 \pm k\pi\right]=(-1)^{k}\sin[(2n+1)t/2]
\]
\subsection{Proof of Theorem \ref{th-LC}}
Regarding the proof that (\ref{eq-LC}) are sufficient condition for (\ref{boundedness}), we refer the reader to \cite[Th. 4.1]{ThBa} for the cases $w\ne w_1$, and to \cite[Eq. (29)]{paperoAMC} for the case $w=w_1$. Here we prove that the bounds in (\ref{eq-LC}) are necessary for having
\[
\sup_{n\sim m}\|V_n^m\|_u< \infty
\]
We obtain the statement  reasoning by contradiction. To this aim, for all $x\in [-1,1]$ we set $t=\arccos x$, so that we can write
\[
u(x)=(1-x)^\gamma (1+x)^\delta=
2^{\gamma+\delta}\left(\sin[t/2]\right)^{2\gamma}\left(\cos[t/2]\right)^{2\delta}
\]
and consequently
\begin{equation}\label{Leb-VP}
\|V_n^m\|_u=\sup_{|x|\le 1}u(x)\sum_{k=1}^n\frac{|\Phi_{n,k}^m(x)|}{u(x_k)}\ge
\sup_{t\ne t_k}\sum_{k=1}^n  \left(\frac{\sin [t/2]}{\sin[t_k/2]}\right)^{2\gamma}\hspace{-.1cm}
\left(\frac{\cos [t/2]}{\cos[t_k/2]}\right)^{2\delta}\hspace{-.2cm}|\Phi_{n,k}^m(\cos t)|
\end{equation}
Moreover, for fixed coprime integers $0<\mu<\nu$, we take the degree--parameters as follows
\begin{equation}\label{deg-mn}
n=2l\nu,\qquad m=2l\mu, \qquad \mbox{with}\qquad l\in\NN,
\end{equation}
and by means of Proposition \ref{prop}, we are going to show that if (\ref{eq-LC}) does not hold, then we get
\begin{equation}\label{Leb-VPsup}
\sup_{n\sim m}
\|V_n^m\|_u \ge \sup_{l\in\NN} \left\{\|V_n^m\|_u : \ n=2l\nu \ \mbox{and}\  m=2l\mu \right\}= +\infty
\end{equation}
Throughout the subsection, wherever we write $\C$, we mean that it is always $0<\C\ne\C(l,t,t_k)$, $k=1,\ldots, n$. Moreover, we write $a\sim b$ to mean that $\C^{-1} a<b<\C \ a$.

Before giving the detailed proof for each Chebyshev case, we remark that in all cases the following estimates hold
\begin{eqnarray}\label{sim-1}
\sin(t_1/2)\sim n^{-1} &\mbox{and}& \cos(t_1/2)\sim 1 \\
\label{sim-2}
\sin(t_n/2)\sim 1\hspace{.4cm} &\mbox{and}& \cos(t_n/2)\sim n^{-1}\\
\label{sim-3}
\sin(t_k/2)\sim t_k \hspace{.2cm}&& \hspace{-1cm}k=1,\ldots,n,
\end{eqnarray}
and $m\sim n$ implies that
\begin{equation}\label{sim-4}
\C\le \sin(t_m/2)\le 1\quad\mbox{as well as}\quad \C\le \cos(t_m/2)\le 1.
\end{equation}
Moreover, by (\ref{diff}) and (\ref{deg-mn}) we easily get
\begin{eqnarray}
\label{diff-1}
{}\hspace{-1cm}\left|\Psi\left(\frac\pi 2, t_1\right)\right|\hspace{-.2cm}&=&\hspace{-.2cm}
|\sin(2l\mu t_1)|\left(
\frac 1{\sin^2\left[\pi/4-t_1/2\right]}+\frac 1{\sin^2\left[\pi/4+t_1/2\right]}
\right)\ge \C,\quad \forall l\in\NN,\\
\label{diff-2}
{}\hspace{-1cm}\left|\Psi\left(\frac\pi 2, t_n\right)\right|\hspace{-.2cm}&=&\hspace{-.2cm}
|\sin(2l\mu t_n)|\left(
\frac 1{\sin^2\left[\pi/4-t_n/2\right]}+\frac 1{\sin^2\left[\pi/4+t_n/2\right]}
\right)\ge \C,\quad \forall l\in\NN
\end{eqnarray}
Furthermore, we note that
\begin{equation}\label{diff-0}
|\Psi(0, t_k)|=2 \frac{|\sin(m t_k)|}{\sin^2(t_k/2)}\quad\mbox{and}\quad
|\Psi(\pi, t_k)|=2 \frac{|\sin(m t_k)|}{\cos^2(t_k/2)},\qquad k=1,\ldots,n,
\end{equation}
as well as we easily get
\begin{equation}\label{diff-m}
\left|\Psi\left(\frac\pi m, t_k\right)\right|\ge 2 |\sin(m t_k)|\quad\mbox{and}\quad
\left|\Psi\left(\pi\hspace{-.1cm}-\hspace{-.1cm}\frac{\pi} m, t_k\right)\right|\ge 2 |\sin(m t_k)|,\hspace{.7cm} k=1,\ldots,n.
\end{equation}
$\bullet$ {\it Case $w=w_1$}

In this case from (\ref{Leb-VP}) and (\ref{fi1-trig}) we deduce that
\begin{equation}\label{LC-1}
\|V_n^m\|_u
\ge \frac \C{n^2} \left(\frac{\sin [t/2]}{\sin[t_k/2]}\right)^{2\gamma}
\left(\frac{\cos [t/2]}{\cos[t_k/2]}\right)^{2\delta}|\cos[nt]|
\left|\Psi(t,t_k)\right|
\end{equation}
holds for all $k=1,\ldots,n$ and any $t\in [0,\pi]$, with $t\ne t_k$.

Hence, let us distinguish the following two cases when (\ref{eq-LC}) does not hold.

Case $(1a)$: Suppose $\gamma>1$ and $\delta \ge 0$.
Then taking $k=1$ and $t=\pi/2$ in (\ref{LC-1}), recalling that $n$ is even, and using (\ref{diff-1}) and (\ref{sim-1}),  we get
\[
\|V_n^m\|_u\ge\frac \C{n^2} \left[
\left(\frac{\sin [t/2]}{\sin[t_k/2]}\right)^{2\gamma}
\left(\frac{\cos [t/2]}{\cos[t_k/2]}\right)^{2\delta}|\cos[nt]|
\left|\Psi(t,t_k)\right|
\right]_{t=\frac\pi 2, \ k=1}\ge \C n^{2\gamma -2}
\]
i.e. (\ref{Leb-VPsup}) holds since $\gamma-1>0$.

Case $(1b)$: suppose $\delta>1$ and $\gamma\ge 0$. Then by (\ref{LC-1}), (\ref{diff-2}) and (\ref{sim-2}), we get
\[
\|V_n^m\|_u\ge\frac \C{n^2} \left[
\left(\frac{\sin [t/2]}{\sin[t_k/2]}\right)^{2\gamma}
\left(\frac{\cos [t/2]}{\cos[t_k/2]}\right)^{2\delta}|\cos[nt]|
\left|\Psi(t,t_k)\right|
\right]_{t=\frac\pi 2, \ k=n}\ge \C n^{2\delta -2}
\]
which implies (\ref{Leb-VPsup}) since $\delta-1>0$.

$\bullet$ {\it Case $w=w_2$}

In this case by (\ref{Leb-VP}) and (\ref{fi2-trig}) we have
\begin{equation}\label{LC-2}
\|V_n^m\|_u
\ge \frac \C m \sup_{t\ne t_k}\sum_{k=1}^n \left(\frac{\sin [t/2]}{\sin[t_k/2]}\right)^{2\gamma}
\left(\frac{\cos [t/2]}{\cos[t_k/2]}\right)^{2\delta}
\frac{|\sin[(n+1)t]|\sin t_k}{(n+1)\sin t}\left|\Psi(t,t_k)\right|
\end{equation}
Hence, let us reason by contradiction and suppose the bounds in (\ref{eq-LC}) do not hold. This can happen in one of the following cases:
\[
\begin{array}{lll}
\mbox{Case}\ (2a):\quad \gamma=0\ \mbox{and}\ \delta\ge 0; & \hspace{.5cm}&\mbox{Case}\ (2d):\quad \delta>3/2\ \mbox{and}\ \gamma> 0;\\ [.1in]
\mbox{Case}\ (2b):\quad \delta=0\ \mbox{and}\ \gamma\ge 0; & \hspace{.5cm}&\mbox{Case}\ (2e):\quad \gamma-\delta<-1\ \mbox{and}\ 0<\gamma,\delta\le 3/2;\\ [.1in]
\mbox{Case}\ (2c):\quad \gamma>3/2\ \mbox{and}\ \delta> 0; & \hspace{.5cm}&\mbox{Case}\ (2f):\quad \gamma-\delta>1\ \mbox{and}\ 0<\gamma,\delta\le 3/2.
\end{array}
\]
In the case $(2a)$, we focus on $t=0$ and by (\ref{LC-2}), (\ref{diff-0}), (\ref{sim-3}), (\ref{sim-4}), we get
\begin{eqnarray*}
\|V_n^m\|_u&\ge&\frac \C m\sup_{t\ne t_k}\left[\sum_{k=1}^n \left(\frac{\cos [t/2]}{\cos[t_k/2]}\right)^{2\delta}
\frac{|\sin[(n+1)t]|}{(n+1)\sin t}\sin t_k\left|\Psi(t,t_k)\right|\right]\\
&\ge&\frac\C m \left[\sum_{k=1}^n \left(\cos [t/2]\right)^{2\delta}
\frac{|\sin[(n+1)t]|}{(n+1)\sin t}\sin t_k\left|\Psi(t,t_k)\right|\right]_{t=0}\\
&=& \frac \C m \sum_{k=1}^n\sin t_k \frac{|\sin(m t_k)|}{\sin^2(t_k/2)}
\ge \frac \C m \sum_{k=1}^m\cos(t_k/2) \frac{|\sin(m t_k)|}{\sin (t_k/2)}\\
&\ge&\frac \C m \cos(t_m/2)\sum_{k=1}^m \frac{|\sin(m t_k)|}{t_k}\ge  \C \left(\frac{n+1}m\right)
\sum_{k=1}^m \frac 1 k\left|\sin\left(k\frac{m\pi}{n+1}\right)\right|
\end{eqnarray*}
Hence, taking into account that $m\sim n$ and that
\begin{equation}\label{series}
\sum_{k=1}^\infty\frac {|\sin(k\alpha)|}k=\infty, \qquad \forall \alpha\in [0,\pi],
\end{equation}
we conclude that $\sup_{m\sim n}\|V_n^m||_u=\infty$ holds in the case $(2a)$.

The same conclusion applies to the case $(2b)$ too, since similarly to the previous case by (\ref{LC-2}), (\ref{diff-0}), (\ref{sim-3}) and (\ref{sim-4}), we have
\begin{eqnarray*}
\|V_n^m\|_u
&\ge&\frac\C m \left[\sum_{k=1}^n \left(\frac{\sin [t/2]}{\sin[t_k/2]}\right)^{2\gamma}
\frac{|\sin[(n+1)t]|}{(n+1)\sin t}\sin t_k\left|\Psi(t,t_k)\right|\right]_{t=\pi}\\
&\ge& \frac \C m \sum_{k=1}^n\sin t_k \frac{|\sin(m t_k)|}{\cos^2(t_k/2)}
\ge \frac \C m \sum_{k=m}^n\sin(t_k/2) \frac{|\sin(m t_k)|}{\cos (t_k/2)}\\
&\ge&\frac \C m \sin(t_m/2)\sum_{k=m}^n \frac{|\sin(m t_k)|}{\cos (t_k/2)}\ge \frac \C m \sum_{k=m}^n \frac{|\sin[m(\pi- t_k)]|}{\sin[(\pi-t_k)/2]}\\
&=&\frac \C m \sum_{h=1}^{n-m+1} \frac{|\sin(m t_h)|}{\sin(t_h/2)}\ge  \C \left(\frac{n+1}m\right)
\sum_{h=1}^{n-m+1} \frac 1 h\left|\sin\left(h\frac{m\pi}{n+1}\right)\right|
\end{eqnarray*}
Regarding the remaining cases, we observe that by (\ref{LC-2}), for all $k=1,\ldots, n$ and any $t\in [0,\pi]$ with $t\ne t_k$, we have
\begin{equation}\label{LC-2bis}
\|V_n^m\|_u
\ge \frac \C m \left(\frac{\sin [t/2]}{\sin[t_k/2]}\right)^{2\gamma}
\left(\frac{\cos [t/2]}{\cos[t_k/2]}\right)^{2\delta}
\frac{|\sin[(n+1)t]|\sin t_k}{(n+1)\sin t}\left|\Psi(t,t_k)\right|.
\end{equation}
Consequently, in the case $(2c)$ the statement (\ref{Leb-VPsup}) follows from (\ref{LC-2bis}), (\ref{sim-1}) and (\ref{diff-1}) which yield
\begin{eqnarray*}
\|V_n^m\|_u&\ge&\frac{\C}m \left[\left(\frac{\sin [t/2]}{\sin[t_k/2]}\right)^{2\gamma-1}
\left(\frac{\cos [t/2]}{\cos[t_k/2]}\right)^{2\delta-1}
\frac{|\sin[(n+1)t]|}{(n+1)}\left|\Psi(t,t_k)\right|\right]_{t=\frac\pi 2,\ k=1}\\
&\ge& \frac \C{n^2}\left(\frac 1{\sin[t_1/2]}\right)^{2\gamma-1}\ge \C n^{2\gamma-3}
\rightarrow \infty \qquad \mbox{as}\qquad n\rightarrow \infty
\end{eqnarray*}
Similarly, in the case $(2d)$, by (\ref{LC-2bis}), (\ref{sim-2}) and (\ref{diff-2}) we have
\begin{eqnarray*}
\|V_n^m\|_u&\ge&\frac \C m \left[\left(\frac{\sin [t/2]}{\sin[t_k/2]}\right)^{2\gamma-1}
\left(\frac{\cos [t/2]}{\cos[t_k/2]}\right)^{2\delta-1}
\frac{|\sin[(n+1)t]|}{(n+1)}\left|\Psi(t,t_k)\right|\right]_{t=\frac\pi 2,\ k=n}\\
&\ge& \frac \C{n^2}\left(\frac 1{\cos[t_n/2]}\right)^{2\delta-1}\ge \C n^{2\delta-3}
\rightarrow \infty \qquad \mbox{as}\qquad n\rightarrow\infty
\end{eqnarray*}
and this implies (\ref{Leb-VPsup}) since $2\delta-3>0$.

In the case $(2e)$, we consider $t=\pi/m$, $k=n$ and deduce the statement from (\ref{LC-2bis}), (\ref{sim-2}), (\ref{sim-3}) and (\ref{diff-m}) as follows
\begin{eqnarray*}
\|V_n^m\|_u&\ge&\frac{\C}m \left[\left(\frac{\sin [t/2]}{\sin[t_k/2]}\right)^{2\gamma-1}
\left(\frac{\cos [t/2]}{\cos[t_k/2]}\right)^{2\delta-1}
\frac{|\sin[(n+1)t]|}{(n+1)}\left|\Psi(t,t_k)\right|\right]_{t=\frac\pi m,\ k=n}\\
&\ge& \frac \C{n^{2\gamma +1}}\left(\sin[t_n/2]\right)^{1-2\gamma}\left(\cos[t_n/2]\right)^{1-2\delta}\ge \frac \C{n^{2\gamma-2\delta + 2}}\rightarrow \infty \quad\mbox{as $n\rightarrow \infty$}
\end{eqnarray*}
Finally, in the case $(2f)$ we take $t=(\pi-\pi/m)$, $k=1$ and similarly to the previous case, by (\ref{LC-2bis}), (\ref{sim-1}), (\ref{sim-3}) and (\ref{diff-m}) we get
\begin{eqnarray*}
\|V_n^m\|_u&\ge&\frac{\C}m \left[\left(\frac{\sin [t/2]}{\sin[t_k/2]}\right)^{2\gamma-1}
\left(\frac{\cos [t/2]}{\cos[t_k/2]}\right)^{2\delta-1}
\frac{|\sin[(n+1)t]|}{(n+1)}\left|\Psi(t,t_k)\right|\right]_{t=\left(\pi-\frac\pi m\right),\ k=1}\\
&\ge&\frac{\C}{n^2} \left(\frac{\cos [\pi/(2m)]}{\sin[t_1/2]}\right)^{2\gamma-1}
\left(\frac{\sin [\pi/(2m)]}{\cos[t_1/2]}\right)^{2\delta-1}\\
&\ge& \frac{\C}{n^{2\delta-2\gamma + 2}}\rightarrow \infty \quad\mbox{as $n\rightarrow \infty$}
\end{eqnarray*}
$\bullet$ {\it Case $w=w_3$}

In this case by (\ref{Leb-VP}) and (\ref{fi3-trig}) we have
\begin{equation}\label{LC-3}
\|V_n^m\|_u
\ge \frac \C m \sup_{t\ne t_k}\left[\sum_{k=1}^n \left(\frac{\sin [t/2]}{\sin[t_k/2]}\right)^{2\gamma}
\left(\frac{\cos [t/2]}{\cos[t_k/2]}\right)^{2\delta}
\frac{|\cos[(2n+1)t/2]|\cos [t_k/2]}{(n+1)\cos [t/2]}\left|\Psi(t,t_k)\right|\right]
\end{equation}
and the bounds in (\ref{eq-LC}) do not hold in one of the following cases:
\begin{itemize}
\item[(3a)] $\gamma>1$ and $\delta\ge 0$ \hspace{.2cm}
$\rightarrow$ Proof similar to the previous case $(1a)$.
\item[(3b)] $\delta=0$ and $\gamma\ge 0$ \hspace{.2cm}
$\rightarrow$ Proof similar to the previous case $(2b)$.
\item[(3c)] $\delta > \frac 32$ and $\gamma\ge 0$\hspace{.2cm}
$\rightarrow$ Proof similar to the previous case $(2d)$.
\item[(3d)] $\gamma-\delta > \frac 12$ and $0\le \gamma\le 1$, $0<\delta\le \frac 32$\hspace{.2cm}
$\rightarrow$ Proof similar to the previous case $(2f)$.
\end{itemize}
$\bullet$ {\it Case $w=w_4$}

In this case by (\ref{Leb-VP}) and (\ref{fi4-trig}) we have
\begin{equation}\label{LC-4}
\|V_n^m\|_u
\ge \frac\C m \sup_{t\ne t_k}\left[\sum_{k=1}^n \left(\frac{\sin [t/2]}{\sin[t_k/2]}\right)^{2\gamma}
\left(\frac{\cos [t/2]}{\cos[t_k/2]}\right)^{2\delta}
\frac{|\sin[(2n+1)t/2]|\sin [t_k/2]}{(2n+1)\sin[t/2]}\left|\Psi(t,t_k)\right|\right]
\end{equation}
and the bounds in (\ref{eq-LC}) do not hold in one of the following cases:
\begin{itemize}
\item[(4a)] $\gamma=0$ and $\delta\ge 0$ \hspace{.2cm}
$\rightarrow$ Proof similar to the previous case $(2a)$.
\item[(4b)] $\gamma>3/2$ and $\delta\ge 0$ \hspace{.2cm}
$\rightarrow$ Proof similar to the previous case $(2c)$.
\item[(4c)] $\delta > 1$ and $\gamma\ge 0$\hspace{.2cm}
$\rightarrow$ Proof similar to the previous case $(1b)$.
\item[(4d)] $\gamma-\delta <- \frac 12$ and $0< \gamma\le \frac 32$, $0\le\delta\le 1$\hspace{.2cm}
$\rightarrow$ Proof similar to the previous case $(2e)$.
\end{itemize}

{\bf Acknowledgments}
This research has been accomplished within Rete ITaliana di Approssimazione (RITA) and partially supported by the GNCS-INdAM funds 2019, project ``Discretizzazione di misure, approssimazione di operatori integrali ed applicazioni''.

{\small{\it{
Donatella Occorsio, {Department of Mathematics, Computer Science and Economics,\\ University of Basilicata,\\  Via dell'Ateneo Lucano 10, 85100 Potenza, Italy }. \\ donatella.occorsio@unibas.it.

{{Woula Themistoclakis, C.N.R. National Research Council of Italy, IAC Institute for Applied Computing ``Mauro Picone'',\\ Via P. Castellino, 111, 80131 Napoli, Italy. \\ {woula.themistoclakis@cnr.it.}}}
  }}}

\end{document}